\documentclass{article}

\usepackage{amssymb,amsmath,amsthm}
\usepackage{mathrsfs}
\usepackage{tikz}
\usetikzlibrary{decorations.pathmorphing}
\usetikzlibrary{patterns,shapes}
\theoremstyle{plain}
\newtheorem{tm}{Theorem}[section]

\theoremstyle{definition}
\newtheorem{defn}[tm]{Definition}

\newtheorem{example}[tm]{Example}
\newtheorem{exercise}[tm]{Exercise}
\newtheorem{exe}[tm]{Exercise}

\makeatletter

\renewcommand\thesubsection{\@Roman\c@section.\@arabic\c@subsection}
\renewcommand{\thetm}{\@Roman\c@section.\arabic{tm}}
\def\l@section#1#2{\ifnum \c@tocdepth >\z@ \addpenalty \@secpenalty \addvspace {1.0em \@plus
 \p@ }\setlength \@tempdima {0em}\begingroup \parindent \z@ \rightskip \@pnumwidth \parfillskip -\@pnumwidth \leavevmode \bfseries \advance \leftskip \@tempdima \hskip -\leftskip #1\nobreak \hfil \nobreak \hb@xt@ \@pnumwidth {\hss #2}\par \endgroup \fi}

\def\@fnsymbol#1{$\dagger$}
\makeatother

\def\C{\mathbb{C}}

\def\R{\mathbb{R}}

\def\Z{\mathbb{Z}}
\def\Id{{\rm id}}
\let\e\varepsilon

\let\o\omega
\let\phi\varphi
\def\id{{\rm id}}
\let\hat\widehat
\def\b{\mathfrak b}

\begin{document}

\author{Yulia Kuznetsova\thanks{This work was partially supported by the French ``Investissements d'Avenir'' program, project ISITE-BFC (contract ANR-15-IDEX-03).}}
\title{Quantum semigroups: A survey}

\maketitle


\abstract{\small Notes of a 8h course given at the University of G\"oteborg during an Erasmus exchange visit, June 11-15, 2018. It is intended for PhD and graduate students familiar with $C^*$-algebras but not specializing in quantum groups. The proofs, if included, are presented in detail, but most facts are not proved. A special emphasis is put onto motivation of the theory and the intuition behind.}

\tableofcontents

\section{Lecture I: Introduction and compact quantum groups}


\subsection{Overview}

Quantum groups in their topological setting are noncommutative analogues of function algebras on topological groups. Every class of quantum (semi)groups corresponds to a classical family of (semi)groups, varying by the degree of continuity of the multipication, absence or presence of compactness, presence of additional structures. In my course, I will define and discuss known families of quantum groups and semigroups and try to systemize the work done beyond the main framework of locally compact quantum groups.
	
Much attention will be paid to examples. Among them will be: partial Hadamard matrices \cite{banica}, quantum families of maps \cite{soltan}, quantum $C^*$-algebras of locally compact semigroups \cite{li,aukh}. The time limits will impose the level of details. Finally, I will speak of certain maps between classes mentioned above: Bohr \cite{soltan-bohr} and weakly periodic \cite{daws} compactifications and duals of semigroups with involution \cite{haar2}. 

I am grateful to Lyudmila Turowska for hospitality during my visit. I thank also Biswarup Das and Piotr So\l tan for comments on these notes.

\subsection{General guiding ideas}

To every locally compact group, one can associate at least two natural group algebras: $C_0(G)$, the algebra of continuous functions vanishing at infinity, and $C^*(G)$, the enveloping $C^*$-algebra of $L^1(G)$. They carry much information about $G$, but do not determine it uniquely. For example, $C_0(\Z)\simeq C_0(\Z^2)$, or even $C_0(G)\simeq C_0(F)$ for every pair of finite groups $F,G$ of the same order.

The group $C^*$-algebras can also be isomorphic for non-isomorphic groups. Let $Q$ be the quaternion group: $Q=\{\pm e,\pm i,\pm j,\pm k\}$ with the identity $e$ and relations $i^2=j^2=k^2=ijk=-e$, the minus acting in the multiplication in the usual way. Next, the dihedral group $D_4$ is the group of symmetries of a square; it has order 8 and includes  four rotations, two reflexions with respect to the axes and two reflexions with respect to the diagonals. These groups are of course not isomorphic because $Q$ has six elements of order 4 and $D_4$ only three. However, one can show that $C^*(Q)\simeq C^*(D_4)$. (This does not follow from cardinality alone, as for example $C^*(\Z_8)$ is commutative and not isomorphic to the two algebras above.)

One could say thus that we need both algebras to characterize $G$, but another option is to consider two structures on any of these algebras: multiplication and comultiplication. Let $G$ be compact. On $C(G)$, there is a map dual to the multiplication: $\Delta: C(G)\to C(G\times G)$,
\begin{equation} \label{Delta}
\Delta(f)(s,t)=f(st).
\end{equation}
If now we have an isomorphism $\phi:C(G)\to C(H)$ which is compatible with comultiplication, then $G$ must be isomorphic to $H$.

The space $C(G\times G)$ is isomorphic to the minimal $C^*$-tensor product $C(G)\otimes C(G)$, so that a simple tensor $f\otimes g$, $f,g\in C(G)$, corresponds to the function $(s,t)\mapsto f(s) g(t)$ on $G\times G$. This tensor product representation is more convenient since it is given in terms of the initial algebra.

$C^*(G)$ has also a natural comultiplication $\Delta:C^*(G)\to C^*(G)\otimes C^*(G)$ which is dual to the pointwise multiplication of functions. To guess the correct formula, let us suppose that $G$ is discrete. For every $s\in G$, let $\delta_s$ be the indicator function of $s$; it belongs to $\ell^1(G)$ and thus to $C^*(G)$ (and acts on $\ell^2(G)$ as a translation operator).  Let us view $\ell^1(G)$ as the linear dual space to $c_0(G)$, and $\ell^1(G\times G)$ as the dual of $c_0(G\times G)\simeq c_0(G)\otimes c_0(G)$. We would like to have $\Delta(\delta_s)(f\otimes g)=\delta_s(fg)$ for all $f,g\in c_0(G)$.

Since $\delta_s(fg) =(fg)(s) = f(s) g(s) = \delta_s(f) \delta_s(g)=(\delta_s\otimes \delta_s)(f\otimes g)$, this implies
$$
\Delta(\delta_s)=\delta_s\otimes\delta_s.
$$ 
By verifying continuity, one checks that in fact, this map extends to a *-homomorphism from $C^*(G)$ to $C^*(G)\otimes C^*(G)$. In the non-discrete case the reasoning requires more details but follows the same ideas.

\subsection{Compact quantum semigroups}

We come to the following definition:

\begin{defn}
A {\it compact quantum semigroup} (CQS) is a pair $(A,\Delta)
$ where $A$ is a unital $C^*$-algebra and the comultiplication $\Delta:A\to A\otimes A$ is a unital *-homomorphism which is coassociative, that is, such that
\begin{equation}\label{coassoc}
(\Id\otimes\Delta)\Delta=(\Delta\otimes\Id)\Delta.
\end{equation}
\end{defn}

To discuss the coassociativity condition, suppose that $A=C(G)$ for a compact group $G$. Let us calculate first the maps $\Id\otimes\Delta$ and $\Delta\otimes\Id$ from $C(G)\otimes C(G)\simeq C(G\times G)$ to $C(G)\otimes C(G)\otimes C(G)\simeq C(G\times G\times G)$. On simple tensors, we have
$$
(\Id\otimes\Delta)(f\otimes g)(s,t,r) = f(s) g(tr) = (f\otimes g)(s,tr),
$$
with $f,g\in C(G)$, $s,t,r\in G$. This implies that
$$
(\Id\otimes\Delta)(F)(s,t,r) = F(s,tr)
$$
for any $F\in C(G\times G)$. Similarly,
$$
(\Delta\otimes\Id)(F)(s,t,r) = F(st,r).
$$
Thus,
$$
(\Id\otimes\Delta)\Delta(f)(s,t,r) = \Delta(f)(s,tr) = f\big( s(tr)\big),
$$
and respectively
$$
(\Delta\otimes\Id)\Delta(f)(s,t,r) = f\big( (st)r\big).
$$
We see now that \eqref{coassoc} holds for every $f\in C(G)$ since the multiplication on $G$ is associative.

Let now, more generally, $A$ be a {\sl commutative\/} CQS. Then, as a unital commutative $C^*$-algebra, $A$ is isomorphic to $C(P)$ on a compact space $P$. For every $s,t\in P$, the map $\phi_{s,t}:A\to\C$, $f\to \Delta(f)(s,t)$ is a *-homomorphism, so that there exists $u_{s,t}\in P$ such that $\phi_{s,t}(f)=f(u_{s,t})$ for every $f\in C(P)$. Let us denote $u_{s,t}$ by $st$, then the formulas above apply and show that this operation $(s,t)\mapsto st$ is associative. It is also not difficult to show that it is jointly continuous. Thus, $P$ is a compact topological semigroup; this explains the term CQS, and suggests that it is more precise to speak of {\it compact quantum topological semigroups (CQTS)}.

\subsection{Compact quantum groups}

To encode the invertibility, one adds the following {\it cancellation rules } to the definition:

\begin{defn}\label{CQG}
A {\it compact quantum group (CQG)} is a CQS such that the sets
\begin{equation}\label{cancel-delta-1}
\{(A \otimes 1) \Delta(A)\} = \{ (a\otimes 1)\Delta(b): a,b\in A\}
\end{equation}
and
$$
\{(1\otimes A) \Delta(A)\} = \{ (1\otimes a)\Delta(b): a,b\in A\}
$$
are dense in $A\otimes A$.
\end{defn}

To understand the idea behind these conditions, let us suppose again that $A=C(P)$ with a compact semigroup $P$. For $f,g\in C(P)$, $s,t\in P$,
$$
\big((f\otimes 1)\Delta(g)\big)(s,t) = f(s) g(st), \quad \big((1\otimes f)\Delta(g)\big)(s,t) = f(t) g(st).
$$
Let us prove that the assumptions of Definition \ref{CQG} are equivalent to the cancellation rules:
\begin{equation}\label{cancel}
st=sr \ \Rightarrow\ t=r, \quad ts=rs \ \Rightarrow\ t=r
\end{equation}
(left and right cancellation). If the left cancellation does not hold, there exist $s,t,r\in P$ such that $st=sr$, but $t\ne r$. By calculations above, the set \eqref{cancel-delta-1} is contained in the kernel of the nonzero functional $\psi:C(P\times P)\to \C$, $\psi(F)=F(s,t)-F(s,r)$, thus it cannot be dense in $C(P\times P)$. 

Suppose now that the left cancellation law in \eqref{cancel} holds. To show that the set \eqref{cancel-delta-1} is dense in $C(P\times P)$, by Stone-Weierstrass theorem it is sufficient to prove that it separates points, as it contains constants and is clearly a *-subalgebra. Now, it is obvious that a pair of points $(s,t)$, $(s',t')$ with $s\ne s'$ is separated by \eqref{cancel-delta-1}: one can set $g\equiv 1$ and pick any $f$ separating $s$ and $s'$. If now $s=s'$, but the points are distinct, then $t\ne t'$. By the cancellation law, $st\ne st'$, so there exists $g\in C(P)$ such that $g(st)\ne g(st')$. It remains now to take $f\equiv 1$.

Most common compact quantum groups belong to the following smaller class:

\begin{defn}\label{CMQG}
A {\it compact matrix quantum group} is a pair $(A, u)$, where $A$ is a unital $C^*$-algebra and $u \in M_n(A)$ is an invertible matrix whose entries generate $A$ in the sense that the *-subalgebra $\cal A$ generated by $(u_{ij})$ is dense in $A$, and such that\\
$\bullet$ the matrix $u^t=(u_{ji})$ is also invertible;\\
$\bullet$ the map $\Delta(u_{ij}) = \sum_{k=1}^n u_{ik}\otimes u_{kj}$ extends to a *-homomorphism from $A$ to $A\otimes A$.
\end{defn}

The idea behind this definition is best seen on the example of a compact algebraic Lie group, such as $SU(n)$, $O(n)$ for example. Such a group $G$ is embedded into $GL(n,\C)$ for some $n$, and there are $n^2$ coordinate functions $(u_{ij})$ on $G$. In particular, $u_{ij}\in C(G)$ for every $i,j$.

On the product of two matrices $g=(g_{ij})$, $h=(h_{ij})$ the values of $u_{ij}$ are as follows:
$$
u_{ij}(gh) = \sum_{k=1}^n g_{ik} h_{kj} = \sum_{k=1}^n u_{ik}(g) u_{kj}(h) = \big(\sum_{k=1}^n u_{ik}\otimes u_{kj}\big)(g,h).
$$
According to the formula \eqref{Delta}, we should have thus
\begin{equation}\label{Delta-uij}
\Delta(u_{ij}) = \sum_{k=1}^n u_{ik}\otimes u_{kj}.
\end{equation}
This formula appears thus in many definitions of quantum groups and semigroups.

If we do not impose the invertibility condition on the generating matrix in Definition \ref{CMQG}, we get $C^*$-bialgebras which are not CQG but still have much in common with matrix groups. They are examples of CQS.

\subsection{Free unitary quantum group $U_n^+$}

Many quantum groups are defined as universal $C^*$-algebras with relations. This is the case, in particular, in the following definition of A.\ Van Daele and S.\ Wang \cite{vd-wang}:

\begin{defn}
The free unitary quantum group $U_n^+$ is the universal unital $C^*$-algebra $A$ generated by the $n^2$ elements $(u_{ij})_{1\le i,j\le n}$ under condition that $u=(u_{ij})$ and $u^t$ are unitary. The comultiplication on $A$ is defined by the formula \eqref{Delta-uij}.
\end{defn}

Let us write in more detail the unitarity conditions $uu^*=u^*u=I_n$, $u^t u^{t*}=u^{t*}u^t=I_n$, where $I_n$ denotes the $n\times n$ identity matrix:
\begin{align*}
\sum_{k=1}^n u_{ik} u^*_{jk} &= \sum_{k=1}^n u^*_{ki}u_{kj} = \delta_{ij}1,\\
\sum_{k=1}^n u_{ki} u^*_{kj} &= \sum_{k=1}^n u^*_{ik}u_{jk} = \delta_{ij}1.\\
\end{align*}
Note that since the generators do not commute, the second line does not follows from the first one.

The construction hidden by the words ``universal $C^*$-algebra'' is as follows. Suppose we are given a set of generators $\{a_1,\dots,a_n\}$ and a finite set $\cal R$ of relations betweeen them. Typically these are polynomial relations as in the example above.  Consider first the free unital *-algebra $\cal A$ generated by the elements $\{a_i\}$, that is, the set of all polynomial expressions in $1, \{a_i,a_i^*\}$. Next, for every $a\in \cal A$ set $\|a\|_*=\sup \|\pi(a)\|$, where supremum is taken over all  *-representations $\pi$ of $\cal A$ on Hilbert spaces such that $\big(\pi(a_i)\big)$ satisfy the same relations as $(a_i)$. The subtlety is that there might be no such representations at all, or that the norm $\|\cdot\|_*$ is infinite for certain $a\in\cal A$. In these cases one says that the universal $C^*$-algebra $C^*(\{a_1,\dots,a_n\}|{\cal R})$ with given relations does not exist. But if it does, it is defined as the completion of the quotient algebra ${\cal A}/\{a\in{\cal A}: \|a\|_*=0\}$ with respect to $\|\cdot\|_*$.

For the case of $U_n^+$, we can define $\pi:\cal A\to\C$ by $\pi(u_{ij})=\delta_{ij}$, then the matrices $V=\big(\pi(u_{ij})\big)$ and $V^t$ are both equal to the identity matrix, thus are unitary; this shows that the  norm on $\cal A$ is well defined. Further, for every relations-preserving $\pi$ we have $\sum_k \pi(u_{ik})^*\pi(u_{ik}) = 1$, thus $\|\pi(u_{ik})\|\le1$ for all $i,k$. This implies that $\|a\|_*<\infty$ for all $a\in \cal A$, and the universal $C^*$-algebra $A$ is well defined. This is a compact matrix quantum group, one of the first constructed.

\subsection{The quantum group $SU_q(2)$}

This was the first ``truly quantum'' group constructed in the $C^*$-algebraic setting. The exact meaning of this word can be explained a bit later.

Let $SU(2)$ be, as usual, the special unitary group, consisting of complex unitary matrices with the determinant 1. Every matrix $u\in SU(2)$ can be written as
$$
u=\begin{pmatrix} a&-\bar c\\ c&\bar a\end{pmatrix},
$$
with $|a|^2+|c|^2=1$.
The algebra of polynomials (of matrix coeffiients) on $SU(2)$ is generated by $a,\bar a,c,\bar c$. 

In the line of $q$-deformation approach of Drinfeld and Jimbo, let now the generating elements $q$-commute: $ac=qca$. Then the matrix
%
\begin{equation}\label{uq}
u=\begin{pmatrix} a&-qc^*\\ c&a^*\end{pmatrix}
\end{equation}
is unitary if the following relations hold:
\begin{align} \label{acq}
a^*a+c^*c=1, \quad &aa^*+q^2 cc^*=1;\\
c^*c=cc^*, \quad ac&=qca, \quad ac^*=qc^*a. \notag
\end{align}
We can consider then the algebra $\cal A$ generated by $a,c$ with the relations above as the algebra of ``non-commuting polynomials'' on the quantum group $SU_q(2)$.

S.\ L.\ Woronowicz \cite{woron-matrix} put this into the $C^*$-framework:
\begin{defn}
For $q\in [-1,1]$, let $A$ be the universal unital $C^*$-algebra $A$ generated by the elements $a,c$ subject to relations \eqref{acq}. Defined as
$$
\Delta(a) = a\otimes a - qc^*\otimes c, \qquad \Delta(c) = c\otimes a+a^*\otimes c,
$$
the map $\Delta$ extends to a *-homomorphism from $A$ to $A\otimes A$. The pair $(A,\Delta)$ is called the quantum group $SU_q(2)$.
\end{defn}

Of course one needs to verify that $A$ exists to justify this definition. The continuity of $\Delta$ follows from the fact that $\Delta(a)$ and $\Delta(c)$ also satisfy relations \eqref{acq}, thus $\|\Delta(p)\|\le\|p\|_*$ for every polynomial $p\in\cal A$.

$SU_q(2)$ is a CQG, and even a CMQG, with the generating unitary matrix is given by \eqref{uq}.

\begin{exercise}
Show that $\Delta(a)$ and $\Delta(c)$ as defined above satisfy the relations \eqref{acq}.
\end{exercise}

\subsection{Quantum permutation group $S_n^+$}

The group $S_n$ of permutations of $n$ points can be realized as a matrix group: it is isomorphic to the group of all matrices $a=(a_{ij})$ in $M_n(\C)$ such that $a_{ij}\in\{0,1\}$ and each column and row has exactly one 1, that is,
$$
\sum_{k=1}^n a_{ik} = \sum_{k=1}^n a_{ki}=1, \ \forall i.
$$
In the $C^*$-setting, the natural analogue of $\{0,1\}$ is the set of all projections: $p=p^2=p^*$. In addition, let the elements $a_{ij}$ be non-commuting. This leads to the following definition:

\begin{defn}
The quantum permutation group $S_n^+$ is the universal unital $C^*$-algebra $A$ generated by the $n^2$ elements $(u_{ij})_{1\le i,j\le n}$ under relations:
\begin{itemize}
\item every $u_{ij}$ is a projection: $u_{ij}=u_{ij}^*=u_{ij}^2$;
\item $\sum_{k=1}^n u_{ik} = \sum_{k=1}^n u_{kj} = 1$ for every $i,j$.
\end{itemize}
The comultiplication on $A$ is defined by the formula \eqref{Delta-uij}.
\end{defn}

In every column and every row, the projections are orthogonal.

\begin{exe}
Show that the relations above imply the orthogonality of $\{u_{ik}:1\le k\le n\}$ and $\{u_{ki}:1\le k\le n\}$ for every fixed $i$.
\end{exe}

\section{Lecture II: Compact quantum semigroups}

\subsection{Typical examples of semigroups}

Before turning to quantum semigroups, let us recall some best known examples of semigroups, with topology or without it.

\begin{itemize}
\item The matrix algebra $M_n(\C)$ is a semigroup, but not a group with matrix multiplication. Its quantum version exists in a purely algebraic setting, but not in the setting of operator algebrbas, so we will not consider it.
\item Maps on finite sets, such as the set of all maps from a finite set $X$ to itself. The quantum family of all maps is well defined, we will discuss it in \S\ \ref{sec-qu-maps}. Another semigroup of similar nature is of partial bijections on a finite set,  its quantum version is defined in \S\ \ref{sec-part-Hadamard}.
\item Subsemigroups of groups, such as $\R_+$ or $\Z_+$. A construction of a quantum semigroup associated to a subsemigroup of a locally compact group is presented in \S\ \ref{sec-subsemigroup}.
\item Topological semigroups not embeddable in a semigroup. An example is given by the segment $[0,1]$ with the operation $(s,t)\mapsto \min(s,t)$. The algebra of continuous functions on such a semigroup is a quantum semigroup, but no quantum-specific examples are known.
\item Semitopological semigroups, that is, semigroups with separately continuous multiplication. This will be discussed in \S\ \ref{sec-semitop}. For the moment, let us just say that this is impossible in the group case: on a (compact) group, separate continuity automatically implies joint continuity.
\end{itemize}

\subsection{Quantum partial permutations}\label{sec-part-Hadamard}

Let $X$ be a set of $n$ elements. Let us consider {\it partial permutations} on $X$: maps $\phi:Y\to X$ defined on a subset $Y\subset X$, which are bijections of $Y$ with its image. This is a semigroup; the composition is well defined since we allow $Y$ to be empty.

To every partial permutation  one can associate a $n\times n$ matrix $a=(a_{ij})$: set $a_{ij}=1$ if $\phi(i)=j$ and $a_{ij}=0$ if $\phi(i)$ is not defined. The set of matrices thus obtained is decribed as all matrices with entries in $\{0,1\}$ having at most one 1 in each column and row; on the contrary of permutations matrices, the columns and rows might not sum up to 1.

This reasoning leads to the following definition given by Banica and Skalski \cite{banica}:
\begin{defn}
The quantum semigroup of partial permutations $\widetilde S_n^+$ is the universal unital $C^*$-algebra $A$ generated by the $n^2$ elements $(u_{ij})_{1\le i,j\le n}$ under relations:
\begin{itemize}
\item every $u_{ij}$ is a projection: $u_{ij}=u_{ij}^*=u_{ij}^2$;
\item $u_{ij}u_{ik}=u_{ji}u_{ki}=0$ for every $i$ if $j\ne k$.
\end{itemize}
The comultiplication on $A$ is defined by the formula \eqref{Delta}.
\end{defn}

The second condition means that the projections $(u_{ij})$ are orthogonal for fixed $i$ and varying $j$, or for fixed $j$ and varying $i$. This implies that
\begin{equation}\label{ule1}
\sum_{k=1}^n u_{ik}\le 1 \text{ and } \sum_{k=1}^n u_{ki} \le 1 \text{ for every } i,
\end{equation}
but maybe not equal to 1 as in the case of the quantum permutation group $S_n^+$.

\subsection{Semigroup $C^*$-algebras}\label{sec-subsemigroup}

There are several ways to associate $C^*$-algebras to semigroups. This question is not so obvious: if $P$ is a sub-semigroup of a discrete group $G$ which generates $G$ (for example, if $G=P\cup P^{-1}$), then the sub-$C^*$-algebra of $C^*(G)$ generated by $\{\delta_s: s\in P\}$ is all of $C^*(G)$.

To obtain semigroup-specific algebras, one considers representations of $P$ by isometries and not by unitary operators. We will consider only one construction, of the reduced semigroup $C^*$-algebra due to X.\ Li \cite{li} (he considered the discrete case). In the paper of Li one can find a discussion of advantages and disadvantages of different approaches; the motivation for his construction was that the algebra of a semigroup as simple as $\Z_+$ should not be ``too big'', and at least it should be nuclear. M. Aukhadiev and the author \cite{aukh} extended this construction to the case of a closed semigroup of a locally compact group.

Let now $G$ be a locally compact group, and let $P\subset G$ be a subsemigroup such that
$$
G=P^{-1}P=\{s^{-1}t: s,t\in P\}.
$$
This imposes certain restrictions on $P$: for any $s,t\in P$, the intersection of the sets $sP$, $tP$ must be nonempty. This is the case, for example, for $P=[0,+\infty)^2$, but this is not true for the free semigroup on two generators $P=\langle a,b\rangle$ since $aP\cap bP=\emptyset$.

Consider the set $L^2(P)=\{f\in L^2(G): {\rm supp} f\subset P\}$. The semigroup $P$ acts on $L^2(P)$ by left translations: for $s,t\in P$ and $f\in L^2(P)$, set
$$
T_s f(t) = f(s^{-1}t).
$$
Now, the algebra $C^*_r(P)$ is defined as the $C^*$-subalgebra generated in $B(L^2(P))$ by the operators $T_s$, $s\in P$ (and their adjoints). One should note that every $T_s$ is an isometry, $T_s^*T_s=\Id$, but not necessarily invertible.

The comultiplication is defined as $\Delta(T_s)=T_s\otimes T_s$, by analogy with the case of group $C^*$-algebras. It is a certain work to show that this map extends to a *-homomorphism of $C^*_r(P)$. In particular, there is an additional restriction on the semigroup to make this possible: the independence of ideals in $P$.

Recall first that a subset $E\subset P$ is a right ideal if $EP\subset E$. In $P=[0,+\infty)$, the ideals are described as $[a,+\infty)$, $a\ge0$. In general, the set $sP$ is a right ideal for every $s\in P$, as well as $s^{-1}P=\{t\in P: st\in P\}$ and the subsets obtained from $P$ by a finite number of such multiplications on the left by $s_j$ or $s_j^{-1}$ with $s_j\in P$. The ideals of this kind and called {\sl constructible}.

One can check that $T_sT_s^*=I_{sP}$ (the operator of multiplication by the indicator function of $sP$), and moreover, for every constructible ideal $X$ the operator $I_X$ is a finite product of some $T_s$ or $T_s^*$, $s\in P$. This implies that $\Delta(I_X)=I_X\otimes I_X$.

Now, the independence of ideals in $P$ means that if for some constructible right ideals $X$, $X_1$,\dots,$X_n$ we have $X = \cup_{j=1}^n X_j$ up to a null set, then $X = X_j$ for some $1 \le j \le n$, also up to a null set. It turns out that this condition is necessary in order for comultiplication to be well defined.

We will check the necessity in the case of two ideals. Suppose that $X = Y\cup Z$ up to a null set, then $I_X=I_Y+I_Z-I_{Y\cap Z}$. We have
\begin{align*}
\Delta(I_X) &= I_{Y\cup Z}\otimes I_{Y\cup Z}, \\
\Delta(I_Y+I_Z-I_{Y\cap Z}) &= I_Y\otimes I_Y + I_Z\otimes I_Z - I_{Y\cap Z}\otimes I_{Y\cap Z}.
\end{align*}
Let us show that an equality is possible only if $I_Y=I_X$ or $I_Z=I_X$. Set $f=I_{Y\setminus Z} \otimes I_{Z\setminus Y}$, then
$$\Delta(I_X) f = f,$$
$$\Delta(I_Y)f = \Delta(I_Z) f = \Delta(I_{Y\cap Z})f = 0,$$
thus $f=0$ in $L^2(P\times P)$.

The following is an example of a semigroup whose ideals are not independent.
\begin{example}
Consider $P=\{0\}\cup[1;1.5]\cup[2;\infty)$ as a subsemigroup of the group $\R$ with respect to usual addition and the usual topology. Such examples are called \emph{perforated semigroups}, since they are obtained from $\R_+$ by deleting some intervals. Further compute the following ideals:
$$X_1=1+P=\{1\}\cup[2;2.5]\cup [3;\infty ),$$
$$X_2=1.5+P=\{1.5\}\cup[2.5;3]\cup [3.5;\infty ),$$
$$X=-1.5+(1+P)=\{1\}\cup\{1.5\}\cup [2;\infty ).$$
We see that $X=X_1\cup X_2$, but no pair of these ideals coincide up to a null set. Hence, the ideals of $P$ are not independent.  
\end{example}

\subsection{Matrix algebras as universal $C^*$-algebras}

We have seen several universal $C^*$-algebras which have no easy description. One can at the same time represent the simplest $C^*$-algebras as universal ones with generating sets and relations.

\begin{example}
The algebra $A=C^*(\langle a: a=a^2=a^*\rangle)$ is isomorphic to $\C$.

To see this, consider the free (non-topological, not necessarily unital) *-algebra $\cal A$ generated by $a$, that is, the set of all polynomials $p(a,a^*)$, $p\in\C[X,Y]$.
Denote ${\cal R} = \{\pi(a)=\pi(a)^2=\pi(a)^* \}$. By definition, the norm $\|\cdot\|_*$ on $\cal A$ is defined by the set
\begin{equation}\label{set-P}
{\cal P}=\{\pi:{\cal A}\to B(H): H\text{ is a Hilbert space, }\pi\text{ is a *-homomorphism preserving }{\cal R}\}.
\end{equation}
For every $\pi\in\cal P$, $\pi(a)$ is thus a projection, and for every polynomial $p=\sum_{k,l} p_{k,l} X^k Y^l\in\C[X,Y]$,
$$
\|\pi\big(p(a,a^*)\big)\| = \big\|\sum_{k,l} p_{k,l}\, \pi(a)\big\| = \big|\sum_{k,l} p_{k,l}\big|
$$
(unless $\pi(a)=0$). It follows that the homomorphism $\rho:\cal A\to \C$, $p(a,a^*)\mapsto \sum_{k,l} p_{k,l}$ is isometric, $|\rho(v)|=\|v\|_*$ for every $v\in \cal A$, so that the closure of $\cal A$ with respect to $\|\cdot\|_*$ is isomorphic to $\C$.
\end{example}

\begin{example}
Let $A$ be the universal unital $C^*$-algebra generated by $n^2$ elements $\{e_{ij}:1\le i,j\le n\}$ subject to relations
\begin{align}\label{rel-eij}
e_{ij}e_{kl} &= \delta_{jk} e_{il}, \quad \forall i,j,k,l\\
{\cal R}: \qquad \qquad e_{ij}^*&=e_{ji}, \quad \quad \forall i,j \notag\\
\sum_{i=1}^n e_{ii}&=1. \notag
\end{align}
We will show that $A$ is isomorphic to $M_n(\C)$.

Denote as above by $\cal A$ the free *-algebra generated by $\{e_{ij}\}$. Let $E_{ij}$ be the matrix units in $M_n(\C)$. The map $e_{ij}\mapsto E_{ij}$ extends clearly to a (surjective) $\cal R$-preserving *-homomorphism $\rho:{\cal A} \to M_n(\C)$, so that $\|\rho(v)\|\le \|v\|_*$ for all $v\in\cal A$. 

For any other relations-preserving $\pi:\cal A\to B(H)$ and every polynomial $v\in\cal A$, $\pi(v) =\sum v_{ij} \pi(e_{ij})$ is a linear combination of $\pi(e_{ij})$ with canonically defined coefficients $(v_{ij})$. It follows that $\pi(v) = \sum v_{ij} \pi(e_{ij})$ if and only if $\rho(v) = \sum v_{ij} E_{ij}$. This defines a map $\tau: M_n(\C)\to B(H)$ such that $\pi=\tau\circ\rho$, and it is immediate to verify that $\tau$ is a *-homomorphism. As every *-homomorphism between $C^*$-algebras, $\tau$ is contractive, and we have then $\|\pi(v)\|\le\|\rho(v)\|$ for every $v\in\cal A$. This proves that $\|v\|_*=\|\rho(v)\|$ and thus $\rho$ extends to an isomorphism from $A$ onto $M_n(\C)$.
\end{example}

\begin{exe}
Write down the algebra $\prod_{l=1}^m M_{d_l}(\C)$ as the universal $C^*$-algebra with generators and relations.
\end{exe}

\subsection{Quantum families of all maps}\label{sec-qu-maps}

Recall that the permutation group $S_n$ is isomorphic to the set of all bijections of a set $X$ of $n$ elements. The set of all maps of $X$ into itself has no inverse anymore, but is a semigroup.

Its quantum analogue can be defined similarly to the case of the groups of quantum permutations $S_n^+$ and partial permutations $\tilde S_n^+$. In the classical case, every map $\sigma: X\to X$ can be encoded by a matrix $(\sigma_{ij})$ with $(0,1)$ entries, so that $\sigma(k) = i$ if and only if $\sigma_{ik}=1$; every matrix such defined has eactly one 1 in every column, but maybe several 1's in certain rows. The composition of maps corresponds to the multiplication of matrices.

Let now $A$ be the universal unital $C^*$-algebra with generators $p_{ij}$, $1\le i,j\le n$, and relations
\begin{align}\label{rel-pij}
p_{ij}=p_{ij}^2=p_{ij}^*,\\
\sum_{i} p_{ik}=1, \quad \forall k. \notag
\end{align}
With the comultiplication $\Delta(p_{ij})=\sum p_{ik}\otimes p_{kj}$ it becomes a compact quantum semigroup.

This is a particular case of a more general construction of P.\ So\l tan \cite{soltan} of {\it quantum families of all maps} between a pair of $C^*$-algebras. The construction of So\l tan is valid for a pair of $C^*$-algebras where the first one is finitely generated and the second one unital finite-dimensional, but it is a semigroup only if the algebras are the same. We will concentrate below on this case.

Return again to the case of a finite space $X$ of $n$ elements. On the $C^*$-level, it is described by the commutative algebra $C(X)\simeq \C^n$. In the spirit of noncommutative geometry, we can consider a finite-dimensional $C^*$-algebra as a non-commutative ``finite quantum space''. The construction below will behave as the space of all maps from this quantum space to itself.

Let now $A$ be a finite-dimensional $C^*$-algebra. It is known that $A\simeq \prod_{l=1}^m M_{d_l}(\C)$ is isomorphic to a direct product of a finite family of full matrix algebras. As we have seen before, every $M_{d_l}(\C)$ can be represented as the universal $C^*$-algebra with generators and relations. For a direct product, this is also possible, by adding relations which make the identities of every $M_{d_l}(\C)$ orthogonal projections in $A$. For example, $\C^n$ is generated by $n$ elements $a_1,\dots,a_n$ under relations
\begin{equation}\label{aj}
a_j=a_j^2=a_j^*, \quad a_ja_j=a_ja_i \ \ \forall i,j, \quad a_1+\dots +a_n=1.
\end{equation}

So let, in general, $A\simeq \prod_{l=1}^m M_{d_l}(\C)$ be generated by $a_1,\dots,a_n$ under a finite set $\cal R$ of relations between them. We can suppose that $\cal R$ implies that $a_j^*a_j\le 1$ and $a_ja_j^*\le 1$ for every $j$ (as it is in the relations \eqref{rel-eij}).

\section{Lecture III: Overview of different structures}

\subsection{Quantum families of all maps (continued)}

\begin{defn}
The quantum space $\mathbb F(A)$ of all maps from $A$ to itself is defined as the universal $C^*$-algebra generated by $n\cdot(d_1^2+\dots d_m^2)$ elements $f^{(k,l)}_{ij}$ with $1\le k\le n$, $1\le l\le m$, $1\le i,j\le d_l$, such that the block matrices $Y_k={\rm diag}(F^{(k,1)},\dots,F^{(k,m)})$, $1\le k\le n$ with blocks $F^{(k,l)}=(f^{(l,k)}_{ij})_{1\le i,j\le d_l}$ satisfy the relations $\cal R$.
\end{defn}
This algebra is well defined due to the presence of norm restricting relations in $\cal R$.

In the case of $A=\C^n$, for example, we have $m=n$ and $d_1=\dots=d_n=1$, so that $Y_k={\rm diag}(f^{(k,1)},\dots,f^{(k,n)})$; the relations \eqref{aj} imply that
$$
f^{(k,l)}= (f^{(k,l)})^2= (f^{(k,l)})^* \ \ \forall k,l; \qquad f^{(1,l)}+\dots+ f^{(n,l)}=1 \ \ \forall l.
$$
Up to a change of notation, these are just the relations \eqref{rel-pij}, so we arrive at the same algebra as defined above. In the case $n=2$ in particular one can show that $\mathbb F(\C^2)$ is isomorphic to the subalgebra in $C([0,1],M_2(\C))$ of all continuous maps $f$ from $[0,1]$ to $M_2(\C)$ such that $f(0)$ and $f(1)$ are diagonal \cite{soltan}.

The comultiplication could have been written in terms of the generators, but it cannot be expressed, for example, by formulas in $Y_k$ only, so it is better to take another way. First, and this is an essential part of the construction, the space $\mathbb F(A)$ comes with a {\it dual evaluation map}
$$
\Phi: A\to A\otimes \mathbb F(A), \quad a_k\mapsto Y_k.
$$
Note that by definition $Y_k\in \prod_{l=1}^m M_{d_l}\big(\mathbb F(A)\big)$, but this space is also isomorphic to
$$
\prod_{l=1}^m M_{d_l}(\C)\otimes \mathbb F(A)\simeq A\otimes \mathbb F(A).
$$
As for the intuitive meaning of this map, look again at the case of a finite set $X$ of $n$ elements with the usual set $F(X)$ of all maps from $X$ to itself. $\Phi$ is dual to the evaluation map $E:X\times F(X)\to X$, $(x,f)\mapsto f(x)$. $Y_k$ should be viewed as the function on $X\times F(X)$ equal to the ``$k$-th indicator function'' of $f(x)$ on a pair $(x,f)$. 

Now we come to the comultiplication. It should of course correspond to the composition of maps: $(f,g)\mapsto g\circ f$, \ $F(X) \times F(X)\to F(X)$, which evaluates as
$$
E(x,g\circ f) = g(f(x))=E(f(x),g)=E(E(x,f),g)= E\circ(E\times \id)(x,f,g).
$$
In the space of functions, this dualizes to
$$
(\id\otimes\Delta)\Phi=(\Phi\otimes\id)\Phi.
$$
It is shown by So\l tan that such a map $\Delta$ exists and is coassociative. With this structure,  $\mathbb F(A)$  is a compact quantum semigroup, but not a quantum group. This is of course expectable since $F(X)$ is not a group already in the classical case.

The definition and construction can be given quite differently in the terms of universal properties of the algebra $\mathbb F(A)$. This is done in the original paper of So\l tan \cite{soltan} in what concerns the definition, and in \cite{skal-soltan} for the construction.

\subsection{Notations}

In further discussion, we will need the following notions.

Tensor products of $C^*$- and von Neumann algebras as supposed to be known. Unless indicated otherwise, $\otimes$ means the minimal $C^*$ tensor product, and $\bar\otimes$ the von Neumann algebraic tensor product.
Speaking of representations of $C^*$-algebras, we always mean *-representations.

If $A$ is a $C^*$-algebra, then there exists a von Neumann algebra $M$, with a canonical embedding $\imath:A\to M$, such that every *-homomorphism $\phi$ of $A$ into a von Neumann algebra $N$ extends to a normal *-homomorphism $\bar\phi:M\to N$ (so that $\bar\phi\circ\imath=\phi$). In particular, this means that every representation of $A$ extends to a normal representation of $M$. This algebra can be constructed as the closure of $A$ in the universal representation. It is also canonically identified with the second dual of $A$, so usually one denotes it as $A^{**}$. 

If $A$ is a $C^*$-algebra, then $M(A)$ denotes its multiplier algebra. It can be defined in several ways; for example, as the algebra of double centralizers $(L,R)$: $L,R$ are bounded maps on $A$ such that $aL(b)=R(a)b$. To every $a\in A$ one associates the pair $(L_a,R_a)$ where $L_a(b)=ab$ and $R_a(b)=ba$ (check that $cL_a(b)=cab=R_a(c)b$). $M(A)$ is a unital $C^*$-algebra which contains $A$ as an essential ideal: for every nontrivial ideal $I\subset M(A)$, the intersection $I\cap A$ is also nontrivial.

For $A=C_0(X)$, $M(A)=C_b(X)$. For $A=K(H)$, $M(A)=B(H)$.

On the other hand, one can define $M(A)$ as the set of multipliers of $A$ in $A^{**}$:
$$
M(A)=\{x\in A^{**}: ax, xa\in A \text{ for any } a\in A\}.
$$

\subsection{Various structures of quantum (semi)groups}

An overview of quantum (semi)group structures is presented on the scheme below, and the letters appearing on it are as follows:   	
\begin{figure*}[tbh]	
	\begin{tikzpicture}
	\draw [fill=black!3] (-4.1,1) -- (-3.1,3) -- (2,3) -- (4.8,1) -- (4.8,-1.5) -- (-0.9,-4.6) -- (-3,-4.6) -- (-4.1,-1.5) -- cycle [dash pattern=on 20pt off 10pt]; 
	\draw [fill=black!3] (-2.1,2) -- (-1,3) -- (2.2,3) -- (5.2,1) -- (5.2,-1.5) -- (1.9,-4.6) -- (0,-4.6) -- (-2.1,-1.5) -- cycle  [dashdotted]; 
\def\lcqsellipse{(1.6,0) ellipse (6.3 and 2.7)}
\def\cqstsellipse{(-3,0) ellipse (4.6 and 2.7)}
	\draw [fill=green!10]\lcqsellipse; 
	\draw [fill=yellow!10] \cqstsellipse; 
	   \begin{scope}
        \clip \lcqsellipse;
        \fill[orange!20] \cqstsellipse;
    \end{scope}
   \def\lcqgellipse{(1.5,0) ellipse (3 and 2)}
	\draw [fill=blue!10] \lcqgellipse; 
	   \begin{scope}
        \clip \lcqgellipse;
        \fill[magenta!20] \cqstsellipse;
    \end{scope}
	\draw \cqstsellipse; 
	\draw \lcqgellipse; 
	\draw (-4.1,1) -- (-3.1,3) -- (2,3) -- (4.8,1) -- (4.8,-1.5) -- (-0.9,-4.6) -- (-3,-4.6) -- (-4.1,-1.5) -- cycle [dash pattern=on 20pt off 10pt]; 
	\draw (-2.1,2) -- (-1,3) -- (2.2,3) -- (5.2,1) -- (5.2,-1.5) -- (1.9,-4.6) -- (0,-4.6) -- (-2.1,-1.5) -- cycle [dash pattern=on 8pt off 5pt]; 
	\draw (1.6,0) ellipse (6.3 and 2.7); 

	\node [inner sep=0pt] at (0,0) {CQG}; 
	\node [inner sep=0pt] at (3,0) {LCQG}; 
	\node [inner sep=0pt] at (-3,0) {CQTS}; 
	\node at (6.5,0) {LCQS};
	\node at (-6,0) {CQSTS};
	\node at (-2,-4) {QEA};
	\node at (1.6,-4) {QSI};
	\node at (0,-5.3) {HVNA};

	\draw [->,bend left,decorate,decoration={snake,amplitude=.15mm,segment length=2mm,post length=2mm}] (3.7,0.8) to 	 node [pos=0.5,above,sloped] {dual} (3.7,-0.8);
	\draw [->,decorate,decoration={snake,amplitude=.15mm,segment length=2mm,post length=2mm}] (-6,-4.5) to [bend left] node [pos=0.5,above,sloped] {\small WAP} (-5.6,-1.2);
	\draw [->,decorate,decoration={snake,amplitude=.15mm,segment length=2mm,post length=1mm}] (6,-0.5) to [bend left] node [above,midway] {Bohr} (0.1,-0.3);
	\draw [->,bend left,decorate,decoration={snake,amplitude=.15mm,segment length=2mm,post length=2mm}] (0.4,-4.2) to node [pos=0.6,above,sloped] {dual} (3,-3);
	\draw (-8.5,-6) rectangle (8.1,3.8);
	\end{tikzpicture}
	\end{figure*}
\medskip
	
{\def\labelenumi{\arabic{enumi}. }
\begin{enumerate}
\item {\bf CQG} is the class of compact quantum groups, introduced by S.\ L.\ Worono\-wicz \cite{woron-cqg}. A CQG is a pair of a unital $C^*$-algebra $A$ and a comultiplication $\Delta$ which maps $A$ to the minimal tensor product $A\otimes A$ of $A$ with itself. Additional density conditions correspond to the cancellation law and in the commutative case imply the existence of the inverse on the underlying group. The classical example of a CQG is the algebra $C(G)$ of continuous functions on a \underline{compact topological group}\nobreak \ $G$.
\item {\bf LCQG} is the second best known class and stands for locally compact quantum groups, as defined by J.\ Kustermans and S.\ Vaes \cite{kust-vaes}. A LCQG is a $C^*$-algebra $A$ with a comultiplication $\Delta:A\to M(A\otimes A)$ (which takes values in the space of multipliers of $A\otimes A$) and a pair of $\Delta$-invariant weights $\phi,\psi$ on $A$ which correspond to the left and right Haar measure. The classical example is given by the algebra $C_0(G)$ of continuous functions vanishing at infinity on a \underline{locally compact group} $G$.
\item {\bf LCQS} means localy compact quantum semigroups and includes just $C^*$-algebras with comultiplication $\Delta:A\to M(A\otimes A)$, without any additional structure. In the commutative case, we get the algebras isomorphic to $C_0(S)$ on a locally compact space $S$, and the comultiplication induces multiplication on $S$ so that it becomes a \underline{topological semigroup}.
\item {\bf CQTS} includes most known quantum semigroups which are not groups. A compact quantum topological semmigroup (CQTS) is a unital LCQS, with the canonical example of the algebra $C(S)$ of continuous functions on a \underline{compact topological semigroup}. The word ``topological'' in the term CQTS is often omitted, however it is essential: if $A\simeq C(S)$ is a commutative CQTS, then the multiplication on $S$ induced by the comultiplication of $A$ is jointly continuous.
\item Weakening conditions on the mutiplication, we get the class {\bf CQSTS} of compact quantum semitopological semigroups defined by M.\ Daws \cite{daws}. The classical example is the algebra $A=C(S)$ where $S$ is a \underline{compact semitopological semigroup}; multiplication on $S$ is just separately continuous, which means that $\Delta(A)$ is not in $C(S\times S)\simeq A\otimes A$ anymore. In the noncommutative case, one can define the space of separately continuous elements $A\overset{sc}{\otimes} A$ in $A^{**}\otimes A^{**}$, where $A^{**}$ is identified with the enveloping von Neumann algebra of $A$. Now, a CQSTS is a unital $C^*$-algebra $A$ with a comultipication $\Delta:A\to A\overset{sc}\otimes A$.
\item A quantum Eberlein algebra ({\bf QEA}) \cite{das-eberlein} is a $C^*$-algebra $A$ with a comultiplication $\Delta:A\to A^{**}\otimes A^{**}$ equipped with a canonical corepresentation on a Hilbert space $H$: as it is often done in the theory of quantum groups, this corepresentation is encoded by its ``generator'' $V\in B(H)\otimes A^{**}$ and is such that the map $\mu \mapsto ({\rm id}\otimes\mu)(V)$ from $A^*$ to $B(H)$ is a homomorphism. It is known that every unital QEA is a CQSTS \cite{das-mroz}.
\item {\bf QSI} is the class of quantum semigroups with involution \cite{haar2}. An involution on a semigroup is a map ${}^*: S\to S$ such that $x^{**}=x$ and $(xy)^*=y^*x^*$ for all $x,y\in S$; if $S$ is a group, then an involution is given by the group inverse. A classical example of a QSI is the algebra $C(S)^{**}$, where $S$ stands for a \underline{compact semigroup with involution} (there is no assumption of compactness in a general QSI however). In general, a QSI is a von Neumann algebra $M$ with a comultiplication $\Delta:M\to M\otimes M$ and a (possibly unbounded) coinvolution $\varkappa$ on $M$, satisfying certain axioms. In the commutative case, $\varkappa$ sends a function to its composition with the involution.
\item Finally, all these classes are embedded in the class {\bf HVNA} of Hopf--von Neumann algebras. A HVNA is a von Neumann algebra $M$ with a comultiplication $\Delta:M\to M\otimes M$; for every $C^*$-algebra $A$ in one the classes defined above (except QSI which is already defined for von Neumann algebras), its enveloping von Neumann algebra $A^{**}$ is a HVNA. 
\end{enumerate}
}

Observing this panorama, one can notice that there is no notion of topological quantum group as such, which is not supposed to be locally compact. To work in this setting, one would need to take a quite different approach, since every commutative $C^*$-algebra is described as the algebra $C(X)$ on a locally compact space $X$.

In conclusion, let us list a few maps between the classes above.

\begin{itemize}
\item The duality is a functor on the category of locally compact quantum groups \cite{kust-vaes}, the second dual of a LCQG $\mathbb G$ being isomorphic to $\mathbb G$.
\item The quantum Bohr compactification \cite{soltan-bohr} sends any locally compact semigroup to a compact quantum group and has the universal property of factoring every morphism to a CQG.
\item The weakly almost periodic (WAP) compactification \cite{daws} sends any locally compact quantum group to a compact quantum semitopological semigroup, having the respective universal property.
\item The duality of quantum semigroups with involution \cite{haar2} sends a QSI to another QSI, which is also a LCQS. The composition of the square of this map with WAP, applied to the algebra $C_0(G)$ on a locally compact group $G$, equals to the enveloping von Neumann algebra of the Bohr compactification.
\end{itemize}

\subsection{Locally compact quantum semigroups}

If we consider a non-compact locally compact group $G$, the algebra $C_0(G)$ is no more unital. This has also influence on the comultiplication. For $f\in C_0(G)$, we have $\Delta(f)(s,s^{-1}t)=f(t)$ for every $s,t\in G$; if $f$ does not vanish in a certain point $t\in G$, then $\Delta(f)$ takes the same value $f(t)\ne0$ on the set $\{(s,s^{-1}t):s\in G\}$ which is homeomorphic to $G$ and not compact, so that $\Delta(f)\notin C_0(G\times G)\simeq C_0(G)\otimes C_0(G)$.

To deal with non-unital algebras, we have therefore to define comultiplication differently. In the example above, every $\Delta(f)$ is a continuous bounded function on $G\times G$; but as it is known, the space $C_b(X)$ of continuous bounded functions on a locally compact space $X$ is isomorphic to the algebra $M(C_0(X))$ of multilpiers of $C_0(X)$. This gives a framework which is valid in the general case: on a $C^*$-algebra $A$, the comultiplication is a map from $A$ to $M(A\otimes A)$.

In order to define coasociativity, we need to speak of the following extensions. Recall that for a $C^*$-algebra $B$, $M(B)$ is canonically embedded into $B^{**}$. A *-homomorphism $\phi:A\to M(B)$ can be thus uniquely extended to a normal *-homomorphism $\tilde \phi:A^{**}\to B^{**}$; $\phi$ is said to be non-degenerate if $\tilde\phi$ is unital.

By definition, $\Delta$ can be thus extended to a unital *-homomorphism $\tilde\Delta:A^{**}\to (A\otimes A)^{**}$, and in particular, we can speak of its restriction to $M(A)$ which is also unital. Moreover, $\tilde\Delta\big(M(A)\big)$ is contained in $M(A\otimes A)$ (this is proved by continuity).

\begin{defn}
A {\it locally compact quantum semigroup} is a pair $(A,\Delta)$, where
\begin{itemize}
\item $A$ is a $C^*$-algebra;
\item $\Delta:A\to M(A\otimes A)$ is a non-degenerate *-homomorphism which is coassociative:
$$
(\id\otimes\tilde\Delta)\Delta = (\tilde\Delta\otimes\id)\Delta.
$$
\end{itemize}
\end{defn}

The definition above applies also to the algebras $C^*(G)$ and $C^*_r(G)$ on a locally compact group $G$.

If $A$ is a {\it commutative\/} LCQS, then, similarly to the compact case, one can show that $A$ is isomorphic to $C_0(P)$ where $P$ is a locally compact topological semigroup. 

\subsection{Weights}

To speak of locally compact quantum groups, it is necessary to say a few words about weights.

\begin{defn}
A {\it weight} $\phi$ on a $C^*$-algebra $A$ is a function $\phi: A^+\to[0,+\infty]$ such that
\begin{itemize}
\item $\phi(a+b)=\phi(a)+\phi(b)$, \ $a,b\in A^+$;
\item $\phi(ta)=t\phi(a)$, $a\in A^+$, $t\in\R_+$. 
\end{itemize}
\end{defn}

An example of a weight is of course the Haar integral on $C_0(G)$ or on $L^\infty(G)$, where $G$ is a locally compact group: for $0\le f\in C_0(G)$, we have $\int_G f\ge0$, but it may be infinite. Still, positive linearity holds.

\begin{example}[The Plancherel weight] Recall that the group von Neumann algebra $VN(G)$, dual to $L^\infty(G)$, is the strong closure of the set of convolution operators $\{L_f:f\in L^1(G)\}$, acting as $L_f(g)=f*g$, $g\in L^2(G)$. For $L=T^*T\in VN(G)^+$,
\begin{equation}\label{phiTT}
\phi(T^*T)=\begin{cases} \|f\|_2^2, & \text{ if } T=L_f \text{ for some } f\in L^2(G)\\
+\infty, & \text{otherwise.}\end{cases}
\end{equation}
Let $\delta$ be the modular function of $G$, and let $f^*(t)=\overline{f(t^{-1})} \delta(t^{-1})$ be the usual involution on $L^1(G)$. If in the formula \eqref{phiTT} we have $f\in C_c(G)$, then $L_f^*L_f=L_g$ with $g=f^* *f$, and
$$
\|f\|^2 = \int_G |f(t)|^2dt = \int_G |f(t^{-1})|^2 \delta(t^{-1})dt = \int_G f^*(t) \overline{f(t^{-1})}dt = g(e).
$$
The formula $\phi(L_g)=g(e)$ is valid for any $g\in C_c(G)\cap VN(G)_+$. On an abelian group, this is equal to $\int_{\hat G} \hat g(u)du$.
\end{example}

\begin{defn}
A weight $\phi$ on a von Neumann algebra $M$ is
\begin{itemize}
\item {\it faithful} if $\phi(a)\ne0$ for every nonzero $a\in M_+$;
\item {\it normal} if the set $\{ a\in M_+: \phi(a)\le t\}$ is $\sigma$-weakly closed in $M$ for every $t\in\R_+$;
\item {\it semi-finite} if the linear span of the set of $\phi$-integrable positive elements
$$
\mathscr{M}^+_\phi=\{a\in M_+: \phi(a)<\infty\}
$$
is $\sigma$-weakly dense in $M$;
\item n.s.f. if $\phi$ is normal, semi-finite, and faithful.
\end{itemize}
\end{defn}

\subsection{Locally compact quantum groups}

There has been several approaches to the notion of a locally compact quantum group (LCQG). The now generally accepted one was given by J.\ Kustermans and S.\ Vaes \cite{kust-vaes} and is based on weights. There exist two versions of this definition: in von Neumann and in $C^*$-algebraic setting. We give here only the first one since it is less technical.

\begin{defn}
A {\it locally compact quantum group} in the setting of von Neumann algebras is a quadruple $(M,\Delta,\phi,\psi)$ where \begin{itemize}
\item $M$ is a von Neumann algebra;
\item $\Delta:M\to M\bar\otimes M$ is a comultiplication (normal unital coassociative *-homomorphism);
\item $\phi$ and $\psi$ are n.s.f. weights on $M$, moreover, $\phi$ is left invariant and $\psi$ is right invariant in the sense defined below.
\end{itemize}
\end{defn}

For a locally compact group $G$, the algebras $L^\infty(G)$ and $VN(G)$ are von Neumann algbraic locally compact quantum groups, dual one to another.

Left invariance means the following:
for all $a\in \mathscr M_\phi^+$ and $\o\in M_*^+$ holds
$$
\phi\big( (\o\otimes\id)(\Delta(a)) \big) = \o(1)\phi(a).
$$
For right invariance, one replaces $\o\otimes\id$ by $\id\otimes\o$.

\medskip
Let us see how this formula works in the case of $M=L^\infty(G)$. Application of $\o\otimes\id$, called slicing, corresponds to integration by one variable: for $f\in L^\infty(G)_+\cap L^1(G)$ and $\o\in L^1(G)=M_*$,
$$
(\o\otimes\id)(\Delta(f)) = \int_G \o(t) \Delta(f)(t,s) dt = \int_G \o(t) f(ts) dt.
$$
Next,
$$
\phi\big( (\o\otimes\id)(\Delta(a)) \big) = \int_G \int_G \o(t) f(ts) dt ds = \int_G f \int_G \o = \o(1) \phi(f),
$$
where the calculation was held of course using the left invariance of the measure.

A $C^*$-algebraic locally compact quantum group is a $C^*$-algebra with a comutiplication and a pair of invariant weights, but the conditions on weights are more complicated and we do not explain them here. In the case of a locally compact group $G$, the algebras $C_0(G)$ and $C^*_r(G)$ are LCQG. However, the full group $C^*$-algebra $C^*(G)$ is not a LCQG since its weights are not faithful. It is more exact therefore to speak of {\it reduced locally compact quantum groups}. The same applies to the von Neumann algebraic version: the enveloping von Neumann algebra $C^*(G)^{**}$ of $C^*(G)$, called the Ernest algebra and denoted also $W^*(G)$, is not a LQCG in this definition, for the same reason. There is no notion of a topological quantum {\it group} which would include the cases $C^*(G)$ and $W^*(G)$.

There are of course examples of LCQG which are not compact and do not correspond to classical groups, but it would be too long to introduce them. 

It is also possible to define a LCQG as a $C^*$-algebra with a comultiplication, one invariant weight and an antipode (antipodes are discussed in Section \ref{sec-QSI}). This is equivalent to the $C^*$-algebraic version of the definition above. 

\subsection{Multiplicative unitaries}

To every LCQG, one can associate a certain unitary operator called a {\it multiplicative unitary} which plays an important role in the theory. The class of {\it manageable multiplicative unitaries} introduced by S.\ L.\ Woronowicz \cite{woron-mult} includes all unitaries coming from LCQG in the sense of Kustermans and Vaes, but it is unknown whether the opposite holds. There are in particular no examples which would fit into the first definition but not to the second.

One can in fact define a LCQG as a $C^*$-algebra generated by a manageable multiplicative unitary. This approach is adopted in several articles \cite{sol-wor,meyer-roy}.

The theory of multiplicative unitaries is motivated by the following fact. On $L^2(G\times G)\simeq L^2(G)\otimes L^2(G)$, define the operator $W$ by
$$
W\xi(s,t) = \xi(s,s^{-1}t),
$$
$\xi \in L^2(G\times G)$.
It is immediate to verify that this operator is isometric. Its adjoint
$$
W^*\xi(s,t) = \xi(s,st)
$$
is also isometric, so that $W$ is unitary.

The interest of this operator is that it generates both algebras $L^\infty(G)$ and $VN(G)$ by slicing: for $\o\in {\cal B}(L^2(G))_*$, there exists a bounded functional $\o\otimes\id$ on ${\cal B}(L^2(G\times G))$ defined by $(\o\otimes \id)(T\otimes R)=\o(T)R$. In particular, for $\o=\o_{\xi,\eta}$ with $\xi,\eta\in L^2(G)$, $\o(T)=\langle T\xi,\eta\rangle$ and
\begin{align*}
\big((\o_{x,y}\otimes\id)(T\otimes R)f \big)(t) &= \langle Tx,y\rangle (R f)(t)
= \int_G (Tx)(s)\overline{y(s)} ds (Rf)(t) \\
&= \int_G (Tx)(s) (Rf)(t) \overline{y(s)} ds =
 \int_G (T\otimes R)(x\otimes f)(s,t) \overline{y(s)} ds.
\end{align*}
It follows that 
$$
(\o_{x,y}\otimes\id)(W)(f)(t) =  \int_G W(x\otimes f)(s,t) \overline{y(s)} ds
=  \int_G x(s) f(s^{-1}t) \overline{y(s)} ds = \big((x\bar y)*f\big)(t),
$$
so that $(\o_{x,y}\otimes\id)(W)= L_{x\bar y}$ with $x\bar y\in L^1(G)$. Varying $x,y$ and taking strong closure, we get exactly the algebra $VN(G)$.

From the other hand,
\begin{align*}
\big((\id\otimes \o_{x,y})(T\otimes R)f \big)(t) &= (Tf)(t) \langle Rx,y\rangle
= (Tf)(t) \int_G (Rx)(s)\overline{y(s)} ds \\
&= \int_G (Tf)(t) (Rx)(s) \overline{y(s)} ds =
 \int_G (T\otimes R)(f\otimes x)(t,s) \overline{y(s)} ds.
\end{align*}
It follows that 
$$
(\o_{x,y}\otimes\id)(W)(f)(t) =  \int_G W(f\otimes x)(t,s) \overline{y(s)} ds
=  \int_G f(t) x(t^{-1}s) \overline{y(s)} ds = f(t) (\check x*\bar y)(t),
$$
so that
$(\o_{x,y}\otimes\id)(W) = M_{\check x*\bar y}$ is the operator of multiplication by the function $g= \check x*\bar y$. It is known that the set $\{\check x*y: x,y\in L^2(G)\}=A(G)$ is a subalgebra of $C_0(G)$, dense in the uniform norm; this implies that its strong closure is $L^\infty(G)$.

Such a unitary $W$ exists for every LCQG $M$. Let $H$ denote the Hilbert space where $M$ is realized (in a canonical way defined by the Haar weight), then we have
$$
M = \{ (\o\otimes\id)(W): \o\in {\cal B}(H)_*\}''.
$$
Moreover, if $\hat M$ is the dual of $M$ (we do not discuss duality here), then
$$
\hat M = \{ (\id\otimes\o)(W): \o\in {\cal B}(H)_*\}''.
$$

Now, the approach of Woronowicz, preceded by Baaj and Skandalis, is to take a unitary with such properties that the two algebras above have structures of bialgebras (comultiplication) and are in duality one to another.

More on multiplicative unitaries can be found in the book \cite{timmermann}.

\section{Lecture IV: Compactifications}

\subsection{Almost periodic functions}

On a locally compact group $G$, there is a distinguished class of {\it almost perioic} functions. If we return to the real line, then a function $f:\R\to \R$ is said to be almost periodic if for any $\e>0$ there exists $L>0$ such that every interval $[a,a+L]$ contains an {\it almost period} of $f$: such $T\in[a,a+L]$ that $|f(x+T)-f(x)|<\e$ for all $x\in\R$. An example of an almost periodic function is $f(x)=e^{i\lambda x}+e^{i\mu x}$ with arbitrary $\lambda,\mu$.

This class of functions defined by H. Bohr \cite{bohr-1925,bohr} was much studied since then in view of different applications. Bochner proved the following criterion \cite{levitan}: a function $f$ is almost periodic if and only if the set of translates $\{{}_xf: x\in\R\}$ is compact in $C_b(\R)$. This can be thus taken as a definition, and as such has sense in a very general case, without even assuming any topology:

\begin{defn}
Let $G$ be a group. A bounded function $f:G\to \C$ is almost periodic if the set of all left translates $\{{}_xf: x\in G\}$ (or, equivalently, of all right translates $\{f_x: x\in G\}$) is relatively compact in $\ell^\infty(G)$ with the uniform norm.
\end{defn}

In fact, even boundedness it not necessary to suppose, it follows automatically.
It was proved by J. von Neumann in 1934 \cite{neumann} that the space of all almost periodic functions on a group $G$ is the uniform closure of the coefficients of all finite-dimensional unitary representations of $G$. It is also an algebra with respect to pointwise multiplication.

Continuity does not follow in general from almost periodicity as deined above \cite{shtern}, but if we suppose a function to be continuous and almost periodic, then it is automatically uniformly continuous. Let $AP(G)$ denote the space of all {\it continuous} almost periodic functions on $G$. With the description above, we see that $AP(G)$ is the uniform closure of the space of coefficients of all finite-dimensional {\it continuous} unitary representations of $G$.

This opened another perspective for this theory. $AP(G)$ is a commutative unital algebra; let $\b G$ be its spectrum. This is a compact space, and it is immediate to see that $G$ is mapped continuously into it: $\imath: G\to \b G$. Moreover, if $f\in C(\b G)\simeq AP(G)$ vanishes on $\imath (G)$, then by definition it is zero; thus, $\imath(G)$ is dense in $\b G$. One can show more: $\b G$ is a compact topological group, and has the following {\it universal property}: if $K$ is a compact group and $\phi:G\to K$ a continuous homomorphism, then there exists a continuous homomorphism $\bar\phi:\b G\to K$ such that $\bar\phi\circ\imath=\phi$.
\begin{equation}\label{bG}
\begin{tikzpicture}
	\node (G) at (0,0) {$G$};
	\node (bG) at (4,0) {$\b G$};
	\node (K) at (2,-2) {$K$};
	\draw [->] (G) edge [above,midway] node {$\imath$} (bG);
	\draw [->] (G) edge [below,midway] node {$\phi$} (K);
	\draw [->] (bG) edge [below,midway] node {$\bar\phi$} (K);
\end{tikzpicture}
\end{equation}

We can thus give a definition which resumes this point of view:

\begin{defn}
Let $G$ be a locally compact group. A compact group $\b G$ is the {\it almost periodic}, or {\it Bohr} compactification of $G$ if there exists a continuous homomorphism $\imath: G\to \b G$ with a dense image, such that for any continuous homomorphism $\phi:G\to K$ to a compact group $K$ there exists a continuous homomorphism $\bar\phi:\b G\to K$ such that the diagram \eqref{bG} is commutative: $\bar\phi\circ\imath=\phi$.
\end{defn}

One proves then that the spectrum of $AP(G)$ has the property above, justifying the existence of the Bohr compactification.

Less relevant to our course, but more straightforward the existence can be proved by considering the direct product of all finite-dimensional unitary irreducible representations of $G$; this is done in Dixmier's book \cite[\S~16]{dixmier}.

\subsection{Corepresentations}

Having a double structure on a $C^*$- or von Neumann bialgebra means that we can study it in two ways: as usual by representations and by corepresentations.

In the case of a compact group $G$, the corepresentations of $A=C(G)$ are dual to group representations. Let $\pi:G\to {\cal B}(H)$ be a finite-dimensional continuous representation of $G$ on a Hilbert space $H$ of dimension $n$. In some fixed basis of $H$, let $\big( \pi_{ij}(g)\big)$ be the matrix of $\pi(g)$ for every $g\in G$. We get $n^2$ matrix coefficients $\pi_{ij}\in C(G)$. The fact that $\pi$ is a homomorphism is expressed by
$$
\pi_{ij}(gh) = \big( \pi(g)\pi(h) \big)_{ij} = \sum_{k=1}^n \pi_{ik}(g) \pi_{kj}(h),
$$
that is, for $\pi_{ij}$ we have the usual formula of comultiplication:
\begin{equation}\label{Deltapi}
\Delta(\pi_{ij}) = \sum_{k=1}^n \pi_{ik}\otimes \pi_{kj}.
\end{equation}
When dealing with infinite-dimensional representations, it is more complicated to deal with matrix coefficients, so we will seek for a coordinate-free formula equivalent to \eqref{Deltapi}.

Clearly, $\pi$ can be represented as a matrix $V=(\pi_{ij}) \in M_n(C(G))$. The identification $M_n(C(G))\simeq M_n(\C)\otimes C(G)$ allows to write is as
$$
V = \sum_{i,j=1}^n E_{ij}\otimes \pi_{ij}
$$
with the usual matrix units $E_{ij}$. Applying maps separately to the tensor factors, we get:
\begin{align}\label{DeltaV}
(\id\otimes\Delta)(V) &= \sum_{i,j=1}^n E_{ij}\otimes \sum_{k=1}^n \pi_{ik}\otimes \pi_{kj}
= \sum_{i,j,k=1}^n E_{ik}E_{kj}\otimes \pi_{ik}\otimes \pi_{kj}\\
&= \sum_{i,j,k,l=1}^n E_{ik}E_{lj}\otimes \pi_{ik}\otimes \pi_{lj}
 = \Big( \sum_{i,k=1}^n E_{ik}\otimes \pi_{ik}\otimes \id\Big) 
 \Big( \sum_{j,l=1}^n E_{lj}\otimes\id\otimes \pi_{lj}\Big).   \notag
\end{align}
The first map is just $V\otimes\id$, and the second is ``$V$ acting on the first and third tensor factor''. There is a special {\it leg numbering notation} for this which appears often in the corepresentation theory. 

Let us define leg numbering notation in general, in a possibly infinite-dimensional case. Let $A$, $B$ be $C^*$-algebras and $V\in A\otimes B$. On $B\otimes B$, define the flip map:
$$
\theta: B\otimes B\to B\otimes B, \qquad \theta(b_1\otimes b_2) = b_2\otimes b_1
$$
(one verifies that this is an isometry). For $V\in A\otimes B$, we define elements $V_{12},V_{13}\in A\otimes B\otimes B$ as
$$
V_{12} = V\otimes \id, \qquad V_{13} = (\id\otimes\theta)(V\otimes \id).
$$
Note that $\theta$ does not belong to $B\otimes B$ but acts on it, and in the second equality we have the value of $\id\otimes\theta$ on $V\otimes\id$ and not the composition of operators.

Using this notation, the identity \eqref{DeltaV} is written as
$$
(\id\otimes\Delta)(V) = V_{12} V_{13}.
$$
This is the general form taken as a definition of a corepresentation, including also the non-unital case:

\begin{defn}
Let $A$ be a LCQS. A corepresentation of $A$ on a Hilbert space $H$ is an operator $V\in {\cal B}(H)\otimes M(A)$ such that
\begin{equation}
(\id\otimes\tilde\Delta)(V)=V_{12} V_{13},
\end{equation}
where $\tilde\Delta:M(A)\to M(A\otimes A)$ is the canonical extension of $\Delta$. If $H$ is finite dimensional, then $V$ is said to be finite dimensional.
\end{defn}

In analogy to group representations, one distinguishes the class of {\it strongly continuous} corepresentations; they are $V$ as above under condition $V\in M({\cal K}(H)\otimes A)$, where ${\cal K}(H)$ is the algebra of compact operators on $H$.

\subsection{Quantum Bohr compactification}

In the quantum case we have no more spaces but we have algebras which are ``noncommutative functions'' on virtual quantum spaces. In this approach, the Bohr compactification is represented by $C(\b G)\simeq AP(G)$. Recall that this is the closure of the space of matrix coefficients of {\it representations of $G$}, that is, of {\it corepresentations of $C_0(G)$}.

Now, let $A$ be LCQS (it need not even be a quantum group). The idea is to define its Bohr compactification as the closure of the space of coefficients of all its unitary finite-dimensional corepresentations. We will get a unital $C^*$-algebra, which happens to be a CQG. This is a result of P.\ So\l tan \cite{soltan-bohr}.

To make things work, one needs to impose certain restrictions on corepresentations. Let $H$ be a finite-dimensional Hilbert space and let $\tau$ be the transposition operator on ${\cal B}(H)$.

\begin{defn}
Let $A$ b a LCQS. A finite-dimensional corepresentation $V$ of $A$ is called {\it admissible} if $V$ and its transpose $V^t:=(\tau\otimes \id)(V)$ are invertible in ${\cal B}(H)\otimes M(A)$.
\end{defn}

If $V$ is admissible then $V$ is similar to a unitary corepresentation, but the converse is in general not true.
The counter-examples however are not LCQG, and it is now stated as the {\it admissibility conjecture} that every  unitary finite-dimensional corepresentation of a LCQG is admissible. In \cite{soltan-bohr,daws-bohr,das-daws-salmi} it was proved that the conjecture holds true for unimodular LCQG, that is, such that the left and right Haar weights coincide.   

Next we need to define matrix coefficients. Recall that for a group representation $\pi:G\to {\cal B}(H)$, they are parametrized by $\xi,\eta\in H$ and are defined as functions $\pi_{\xi,\eta}:G\to\C$ such that $\pi_{\xi,\eta}(g)=\langle\pi(g)\xi,\eta\rangle$. If we denote by $\o_{\xi,\eta}$ the linear functional on ${\cal B}(H)$, $\o_{\xi,\eta}(T)=\langle T\xi,\eta\rangle$, then $\pi_{\xi,\eta}=\o_{\xi,\eta}\circ \pi$. 

Now let $G$ be compact, $H$ finite-dimensional, and let us look at the corresponding corepresentation of $C(G)$:
$$
V=\sum_{i,j} E_{ij}\otimes \pi_{ij}\in {\cal B}(H)\otimes C(G).
$$
We can write $\pi_{\xi,\eta}$ as the function
$$
g\mapsto \langle \pi(g)\xi,\eta\rangle = \sum_{i,j} \pi_{ij}(g) \xi_j\bar \eta_i 
= \sum_{i,j} \o_{\xi,\eta}(E_{ij}) \pi_{ij}(g) = \o_{\xi,\eta}\Big(\sum_{i,j} E_{ij}\otimes \pi_{ij}\Big)(g),
$$
what leads to the following coordinate-free form:
$$
\pi_{\xi,\eta} = (\o_{\xi,\eta}\otimes \id)(V).
$$

In the finite-dimensional case, every functional $\o\in {\cal B}(H)_*$ is a linear combination of vector functionals $\o_{\xi,\eta}$. But even in the infinite-dimensional case, one prefers to consider all normal functionals, as this is more natural to define and enlarges the class of coefficients only by uniformly converging series of vector coefficients. 

\begin{defn}
Let $A$ be a LCQS and let $V\in {\cal B}(H)\otimes M(A)$ be its corepresentation. A {\it matrix coefficient} of $V$ is
$$
(\o\otimes \id)(V)\in M(A)
$$
with any $\o\in {\cal B}(H)_*$.
\end{defn}

The following is \cite[Proposition 2.13]{soltan-bohr}:
\begin{tm}
Let $\mathbb G=(A,\Delta)$ be a LCQS. Let $\mathbb {AP}(A)$ be the closure in $M(A)$ of the space of the matrix elements of all admissible finite-dimensional corepresentations of $A$. Then $\mathbb {AP}(A)$ is a unital $C^*$-subalgebra of $M(A)$; $\Delta$ maps $\mathbb {AP}(A)$into $\mathbb {AP}(A)\otimes \mathbb {AP}(A)$, and with this comultiplication $\mathbb {AP}(A)$ is a compact quantum group. It is denoted by $\b \mathbb G$ and is called the {\it Bohr compactification} of $\mathbb G$.
\end{tm}

Thus, the Bohr compactification maps the class of LQCS into the class of CQG.

$\b \mathbb G$ has a universality property similar to the classical one. To formulate it correctly, we need to recall that quantum groups correspond to functional algebras on groups, and this ``reverses the arrows'': if $\phi:G\to H$ is a continuous map, then the dual map $\phi^*: C(H)\to C(G)$, $f\mapsto f\circ\phi$ acts in the opposite direction.

For maps we allow the following class, compatible with multiplication and comultiplication:
\begin{defn}
A morphism of LCQS $(A,\Delta_A)$ and $(B,\Delta_B)$ is a non-degenerate *-homomorphism $\phi:A\to M(B)$ such that
$$
(\tilde\phi\otimes\tilde\phi)\circ\Delta_A=\tilde\Delta_B\circ\phi,
$$
where $\tilde\Delta_B$ is the canonical extension of $\Delta_B$ to $M(B)$ and $\tilde\phi$ is the canonical extension of $\phi$ to $M(A)$ (the restriction of its normal extension to $A^{**}$).
\end{defn}

Every CQG is a LCQS, so that the following statement makes sense:
\begin{tm}
Let $(A,\Delta)$ be a LCQS, and let $(B,\Delta_B)$ be a CQG. If\/ $\phi:B\to A$ is a morphism of LCQS, then there exists a morphism of LCQS\/ $\bar\phi: B\to \mathbb {AP}(A)$ such that $\imath\circ\bar\phi=\phi$.
\begin{equation}\label{qbG}
\begin{tikzpicture}
	\node (A) at (0,0) {$A$};
	\node (APA) at (4,0) {$\mathbb{AP}(A)$};
	\node (B) at (2,-2) {$B$};
	\draw [->] (APA) edge [above,midway] node {$\imath$} (A);
	\draw [->] (B) edge [below,midway] node {$\phi$} (A);
	\draw [->] (B) edge [below,midway] node {$\bar\phi$} (APA);
\end{tikzpicture}
\end{equation}
\end{tm}

\subsection{Semitopological semigroups}\label{sec-semitop}

There are natural cases of semigroups when the multiplication is not jointly continuous. For example, let $\cal U(H)$ be the group of unitary operators on a Hilbert space $H$; if $H$ is infinite-dimensional, then the multiplication is not jointly continuous in weak operator topology. Still, it is separately continuous.

One may say that this example is not too significant since one can consider $\cal U(H)$ in the strong operator topology, and it will become a topological group with a jointly continuous multipilcation. Thus the main motivation for the study of STS was another example.

W. Eberlein \cite{eberlein}, in search for ergodic theorems, initiated study of the following class of functions. Let $G$ be a locally compact group. A continuous function $f$ on $G$ is said to be {\it weakly almost periodic} (w.a.p.) if the set of its left translates $L_xf$, $x\in G$, is relatively compact in $C(G)$ with respect to the weak topology. Clearly, every almost periodic function is w.a.p. But less trivially, every function $f\in C_0(G)$ and every positive definite function is w.a.p., and every w.a.p.\ function is uniformly continuous.

Let $WAP(G)$ denote the space of all w.a.p.\ functions on $G$ (it is closed in the uniform norm). This is a commutative unital Banach algebra. Let $G^{WAP}$ be its spectrum. It turns out that $G^{WAP}$ has a natural structure of a compact semitopological semigroup; this can be shown by viewing $G^{WAP}$ as the weak operator closure of $G$ in the space ${\cal B}(C_b(G))$ of bounded operators on $C_b(G)$. This fact attracted attention to the class of semitopological semigroups.

Considered with the canonical morphism $\imath:G\to G^{WAP}$ with dense image, $G^{WAP}$ has the following universality property: for every continuous homomorphism $\phi: G\to P$ to a compact {\it semitopological semigroup} $P$ there exists a continuous homomorphism $\bar\phi: G^{WAP}\to H$ such that $\bar\phi\circ\imath=\phi$. 

An important property of weakly almost periodic functions (which might not vanish at infinity) is the existence of an invariant mean \cite{glicksberg}: a linear functional $M:WAP(G)\to C$, such that
$$
M(1)=1, \quad M(f)\ge0 \text{ if }f\ge0,
$$
and for every $f\in WAP(G)$
$$
M(f) = M({}_xf)=M(f_x), \quad \forall x\in G.
$$
It is constructed as the Haar integral over a certain compact subgroup of $G^{WAP}$.

\section{Lecture V: Additional weak structures}

\subsection{Hopf--von Neumann algebras}

This is the most general structure among $C^*$- or von Neumann bialgebras:
\begin{defn}
A Hopf--von Neumann algebra (HVNA) is a von Neumann algebra $M$ with a comultiplication $\Delta:\to M\bar\otimes M$ which is a normal *-homomorphism and is coassociative: $(\Delta\otimes\id)\Delta=(\id\otimes\Delta)\Delta$.
\end{defn}

All algebras considered in this course either are already HVNA, or are $C^*$-algebras such that their comultiplication extends to $A^{**}$, which becomes a HVNA.

For every HVNA $M$, the predual space $M_*$ is a Banach algebra with the multiplication generated by $\Delta$: for $\mu,\nu\in M_*$ and $a\in M$, set
$$
(\mu*\nu)(a) = (\mu\otimes\nu)\big(\Delta(a)\big).
$$
In particular, for $M=L^\infty(G)$ and $M_*=L^1(G)$, this is the usual convolution of functions.

Moreover, $M$ and $M_*$ are modules one over another: for $a\in M$, $\mu\in M_*$ we can define $a.\mu,\ \mu.a\in M_*$ by
\begin{equation}\label{a.mu}
(a.\mu)(b)=\mu(ba) \text{ and } (\mu.a)(b) = \mu(ab), \quad b\in M;
\end{equation}
and define $\mu*a$, $a*\mu\in M$ by
\begin{equation}\label{a*mu}
\nu(\mu*a) = (\nu*\mu)(a), \quad \nu(a*\mu)=(\mu*\nu)(a), \quad \nu\in M_*.
\end{equation}

\begin{exe}
Show that for $a\in L^\infty(G)$, $\mu\in L^1(G)$
$$
(a.\mu)(t) = (\mu.a)(t) = a(t)\mu(t),
$$
$$
(\mu*a)(t) = \int_G a(s) \mu(t^{-1}s) ds, \qquad (a*\mu)(t) = \int_G a(st) \mu(s) ds, 
$$
so that the first action is just pointwise multiplication (by $a$), and the second one is the convolution with $\check\mu(t)=\mu(t^{-1})$ on the right and with $\mu(t^{-1})\delta(t^{-1})$ on the left respectively (where $\delta$ is the modular functtion of $G$).
\end{exe}

\subsection{Compact quantum semitopological semigroups}

We aim now to formulate separate continuity in terms of comultiplications, following M. Daws \cite{daws}. Let $P$ be a compact semitopological semigroup and let $f\in C(P)$. We can still define a function $\Delta(f)$ on $P\times P$ as $\Delta(f)(s,t)=f(st)$; it might be discontinuous but it is at least bounded. For fixed $s$, the functions
$$
L_sf = \Delta(f)(s,\cdot) \text{ and } R_sf=\Delta(f)(\cdot,s)
$$
are continuous on $P$. Moreover, one can show \cite{runde} that for any Radon measure $\mu$ on $P$,
$$
t\mapsto \int_P f(st)d\mu(s), \quad t\mapsto \int_P f(ts)d\mu(s)
$$
are continuous functions on $P$.

Let ${\cal M}(P)$ denote the space of all Radon measures on $P$. A natural space where $\Delta$ takes its values would be $C(P)^{**}\otimes C(P)^{**}$ which is the dual to ${\cal M}(P)\otimes {\cal M}(P)$. One can show that the condition
$$
(\mu\otimes \id)(F), \quad (\id\otimes\mu)(F) \in C(P) \qquad \forall \mu\in {\cal M}(P)
$$
implies exactly that $F$ is separately continuous.

This allows to give a non-commutative definition of separate continuity. Recall that every $C^*$-algebra $A$ can be viewed as a subset of its second dual $A^{**}$, identified with its enveloping von Neumann algebra. For a bounded linear map $\phi$ on $A$, let $\tilde\phi$ denote its normal extension to $A^{**}$ (if it exists).

\begin{defn}
A {\it compact semitopological quantum semigroup} is a pair $\mathbb S = (A, \Delta)$ where
\begin{itemize}
\item $A$ is a unital $C^*$-algebra, considered as a norm closed $C^*$-subalgebra of $A^{**}$;
\item $\Delta : A \to A^{**}\bar\otimes A^{**}$ is a unital *-homomorphism satisfying $(\tilde\Delta \otimes \id)\circ\Delta = (\id\otimes \tilde\Delta)\circ \Delta$;
\item for any $a\in A$ and any $\omega\in A^*$,
$$
(\tilde\o\otimes \id)(\Delta(a)) \in A; \quad (\id \otimes \tilde\o)(\Delta(a)) \in A.
$$
\end{itemize}
\end{defn}

\begin{defn}\label{weak-cancel}
Say that a CQSTS $A$ satisfies the {\it weak cancellation laws} if for every state $\o\in A^*$, the linear span of each of the sets
$$
\{ a*(\o.b): a,b\in A\}, \qquad \{ (\o.b)*a: a,b\in A\}
$$
(in notations \eqref{a.mu} and \eqref{a*mu}) is norm dense in $A$.
\end{defn}

\begin{exe}
Show that a CQG satisfies weak cancellation laws.
\end{exe}

It was proved by B.\ Das and C.\ Mrozinski \cite{das-mroz} that if $A$ is a CQSTS which satisfies the weak cancellation laws then $A$ is a compact quantum group.

In fact, one can show that for a compact semitopological semigroup $P$, weak cancellation laws for $C(P)$ hold if and only if $P$ is a group; so that the main point in the theorem of Das and Mrozinski is the joint continuity of multiplication. It is thus a non-commutative generalization of the well-known Ellis continuity theorem: a compact semitopological semigroup which is algebraically a group is a topoloical group.

\subsection{Quantum weakly almost periodic compactification}

On the contrary to Bohr compactification, the w.a.p.\ compactification is not based on the coefficients of representations, since in general even the coefficients of all unitary representations of $G$ might not cover the whole of $WAP(G)$. This suggests another construction in the quantum case also, realized by M. Daws \cite{daws}. As the main property is taken the weak periodicity itself.

Let $G$ be a (classical) locally compact group. For $f\in L^\infty(G)$ and $\mu\in L^1(G)$, define $f_\mu\in L^\infty(G)$ by
$$
\langle \nu,f_\mu\rangle = \langle \mu*\nu,f\rangle = (\mu\otimes\nu)(\Delta(f)),
$$
so that $f_\mu = (\mu\otimes \id)\big(\Delta(f)\big)$. It was proved by \"Ulger \cite{uelger} that $f\in L^\infty(G)$ is w.a.p.\ if and only if the orbit map $L^1(G)\to L^\infty(G)$, $\mu\mapsto f_\mu$ is a weakly compact operator. This has sense in the general case:

\begin{defn}
Let $(M,\Delta)$ be a HVNA. The space of weakly almost periodic elements is defined as
$$
{\rm wap}(M) = \{ x\in M: \text{ the map } M_*\to M, \ \mu\mapsto (\mu\otimes\id)(\Delta(x)) \text{ is weakly compact} \}.
$$
\end{defn}

It can be also defined by slicing in a bigger space, as proved by Daws \cite{daws}. Define the space of ``separately continuous functions''
$$
SC(M\times M) = \{ u\in M^{**}\bar\otimes M^{**}: (\tilde\mu\otimes \id)(u) \in M; \quad (\id \otimes \tilde\mu)(u) \in M, \ \forall \mu\in M^*\},
$$
where $\tilde\mu$ stands for the normal extension of $\mu$ to $M^{**}$. In general, this is not a *-subalgebra of $M^{**}\bar\otimes M^{**}$.
It turns out that $x\in {\rm wap}(M)$ if and only if $\Delta(x)$ belongs to the space of
\[\label{MscM}
M\overset{\rm sc}\otimes M = \{ u\in M\otimes M: u,uu^*,u^*u\in SC(M\times M)\}.
\]

Now we can define the w.a.p.\ compactification:
\begin{defn}\label{def-wap}
Let $(M,\Delta)$ be a HVNA. The space ${\rm wap}(M)$ is a $C^*$-algebra, $\Delta$ maps it into ${\rm wap}(M)\otimes{\rm wap}(M)$, and the pair $({\rm wap}(M),\Delta|_{{\rm wap}(M)})$ is a compact quantum semitopological semigroup. It is called the {\it quantum weakly almost periodical (w.a.p.) compactification} of $M$. 
\end{defn}

It has the expected universality property. To formulate it, we need to define morphisms appearing the statement. Let $(M,\Delta_M)$ be a HVNA, and let $(A,\Delta_A)$ be a CQSTS. A {\it morphism} from $A$ to $M$ is a *-homomorphism $\phi:A\to M$ such that $(\tilde\phi\otimes\tilde\phi) \Delta_A = \Delta_M \phi$, where $\tilde\phi:A^{**}\to M$ is the normal extension of $\phi$.

\begin{tm}
Let $(M,\Delta)$ be a HVNA, and let $(A,\Delta_A)$ be a CQSTS. If\/ $\phi:A\to M$ is a morphism, then there exists a morphism\/ $\bar\phi: A\to {\rm wap}(M)$ such that $\imath\circ\bar\phi=\phi$.
$$
\begin{tikzpicture}
	\node (M) at (0,0) {$M$};
	\node (WAP) at (4,0) {${\rm wap}(M)$};
	\node (A) at (2,-2) {$A$};
	\draw [->] (WAP) edge [above,midway] node {$\imath$} (M);
	\draw [->] (A) edge [below,midway] node {$\phi$} (M);
	\draw [->] (A) edge [below,midway] node {$\bar\phi$} (WAP);
\end{tikzpicture}
$$
\end{tm}

In comparison with the Bohr compactification one can note that for discrete quantum groups, the space of almost periodic elements $\mathbb{AP}(A)$ can be defined by a condition similar to formula \eqref{MscM}: $\mathbb{AP}(A)$ is the set of $x\in M(A)$ such that $\Delta(x)\in M(A)\otimes M(A)$ \cite{soltan-discrete}.

\subsection{Quantum Eberlein algebras}

The motivation to consider this class of algebras is to regain the possibility to work with representations in Hilbert spaces, which is not always possible with semitopological semigroups \cite[Theorem 4.7]{megrel}. In fact, even on a noncompact abelian group $G$ the space $B(G)$ of the coefficients of all continuous unitary representations of $G$ is not dense in $WAP(G)$ (in the uniform norm), which means that this class of representations does not separate point of $G^{WAP}$.

One can define an intermediate compactification $G^{\cal E}$, named after Eberlein, as the spectrum of the uniform closure of $B(G)$. We will not discuss this in detail but look rather at the class of semigroups which includes such compactifications.

Let us call for a moment an {\it Eberlein semigroup} a semigroup $P$ for which there exists a homomorphism and homeomorphism $\phi$ from $P$ to the space $B(H)$ of operators on a Hilbert space $H$, considered in the weak operator topology, with every $\phi(s)$ being a contraction. The set of all contractions in $B(H)$ is itself an example of an Eberlein semigroup; another one is given by the mentioned Eberlein compactification of a locally compact group.

The quantum formulation of such an object takes the following form: 
\begin{defn} A quantum Eberlein algebra is a quadruple $\mathbb S = (A, \Delta, V, H)$ where:
\begin{enumerate}
\item $A$ is a $C^*$-algebra, $\Delta : A \to A^{**}\otimes A^{**}$ is a non-degenerate coassociative *-homomorphism;
\item $V \in B(H)\otimes A^{**}$ is a contraction;
\item $V$ and $\Delta$ are connected in the sense that $(\id\otimes\tilde\Delta)(V) = V_{12}V_{13}$, where $\tilde\Delta$ is the normal extension of $\Delta$ to $A^{**}$;
\item $A_V := \{(\o\otimes\id)(V ) : \o\in B(H)_*\}$ is a subset of $A \subset A^{**}$, and is norm dense in $A$.
\end{enumerate}
\end{defn}
Recall that if $A$ is not unital, then by $\Delta$ being non-degenerate we mean that the normal extension
$\tilde\Delta : A^{**}\to A^{**}\otimes A^{**}$ is unital. Condition 3 says, in other words, that $V$ is a corepresentation of $A$.

A commutative unital QEA $A$ is isomorphic to $C(P)$ on a compact semigroup $P$. Conditions 2 and 3 imply the existence of a homomorphism $\phi: P\to B(H)$ with the range in the set of contractions. In condition 4, the inclusion $A_V\subset A$ corresponds to continuity of the coefficients of the representation, and density means that the representation $\phi$ separates points of $P$.

It is known \cite{das-eberlein} that every unital QEA is a CQSTS. Moreover, by using convex analysis, Das and Daws show \cite[Theorem 4.6]{das-eberlein} that every unital QEA has an invariant mean -- a result that does not exist (for the time being) for the class of CQSTS, where it would be natural to expect it as compared to the classical case.

In the noncommutative case, an invariant mean on a CQSTS $A$ is defined as a state $h\in A^*$ such that for every $a\in A$,
$$
(\id\otimes \tilde h)\big(\Delta(a)\big) = h(a) 1 = (\tilde h\otimes\id)\big(\Delta(a)\big).
$$
In notations \eqref{a.mu} with $M=A^{**}$, this is equivalent to
$$
h(a)\mu(1) = h(a*\mu)=h(\mu*a) \qquad \text{ for all }\mu\in A^*,\ a\in A.
$$

\begin{exe}
Show the equivalence of these two definitions.
\end{exe}

It makes difference whether an invariant mean is faithful ($h(a)=0$, $a\ge0$ implies $a=0$) or not. The invariant mean on $WAP(G)$, $G$ being a locally compact group, is not faithful as it vanishes on $C_0(G)$. It is known \cite{mukher} that a locally compact semitopological semigroup with a faithful invariant mean and a neutral element is a locally compact group.

Das and Mrozinski \cite{das-mroz} have proved a quantum version of this theorem:  they showed that a CQSTS admitting a {\bf faithful} invariant mean and a bounded counit (a homomorphism $\e:A\to\C$ such that $(\tilde\e\otimes\id)\Delta=\id=(\id\otimes\tilde\e)\Delta$) is necessarily a compact quantum group.

\subsection{Quantum semigroups with involution}\label{sec-QSI}

The motivation for considering the following class was a search for a duality theory not based on the Haar weight, as opposed to the duality of LCQG. The appearing class is however interesting in itself and joins the family of topological quantum semigroups.

Let for the moment $P$ be a semigroup. A map $*:P\to P$ is called an {\it involution} if
\begin{itemize}
\item $x^{**}=x$ for all $x\in P$;
\item $(xy)^*=y^*x^*$ for all $x,y\in P$.
\end{itemize}
$P$ equipped with such a map is called a {\it semigroup with involution}. If in addition $P$ is a topological semigroup, we will suppose that its involution is continuous.

An example of involution is of course the group inverse if $P$ is a group. Other examples include:
\begin{itemize}
\item the identity map on any abelian semigroup;
\item the adjoint map on the unit ball of $B(H)$;
\item the weakly almost periodic compactification $G^{WAP}$ of a locally compact group $G$; the natural involution on $G^{WAP}$ extends by continuity the inverse of $G$. 
\end{itemize}

An involution on $P$ generates the dual map $S$ on $C(P)$ which we will call coinvolution:
$$
Sf(t)=f(t^*).
$$
If $P$ is a group, then $S$ is usually called an antipode. It plays an important role in the theory of Hopf algebras, topological or not, even if in the topological case it does not appear in the definitions. On every LCQG, for example, there exists an antipode (with certain properties connecting it to the comultiplication and the weights). And in fact, this is the antipode which distinguished ``truly quantum'' groups from classical ones: the class of LCQG with everywhere defined and bounded antipode includes not so many objects apart from the group algebras. All $q$-deformations of Lie groups, in their turn, have unbounded antipodes.

\begin{exe}
On the generators of $SU_q(2)$, set
$$
S(a)=a^*, \quad S(a^*)=a,\quad S(c)=-qc, \quad S(c^*)=-q^{-1}c^*.
$$
Show that $S$ extended as an anti-homomorphism to the algebra generated by $a,a^*,c,c^*$ is not bounded in the $\|\cdot\|_*$-norm. (Hint: $\|c\|_*=\|c^*\|_*\le1$.)
\end{exe}

The definition below allows thus for an unbounded coinvolution.

\begin{defn}\label{def-coinvolution}
Let $M$ be a HVNA. A linear map $S: D(S)\subset M\to M$ is called a {\it proper coinvolution} if it satisfies the following conditions:
\begin{enumerate}
\item $D(S)$ is $\sigma$-weakly dense in $M$;
\item $D(S)$ is closed under multiplication and $S: D(S)\to M$ is an anti-homomorphism;
\item $(*S)(D(S))\subset D(S)$ and $(*S)^2=\Id_{D(S)}$;
\item if $\mu,\nu\in M_*$ are such that $\mu\circ S$ and $\nu\circ S$ extend to normal functionals on $M$, then for all $x\in D(S)$ holds $(\nu\circ S\otimes\mu\circ S)\big( \Delta (x)\big) = (\mu\otimes\nu) \big( \Delta S(x)\big)$.
\end{enumerate}
A HVNA equipped with a proper coinvolution is called a {\it quantum semigroup with involution}.
\end{defn}

The last formula is a replacement of the identity $\theta(S\otimes S)\Delta=\Delta S$, where $\theta$ is the flip map, which may make no sense a priori. The antipode of a locally compact quantum group satisfies these conditions, see \cite[Lemma 5.25]{kust-vaes}.

\begin{example}\label{ex-semigroup}
Let $P$ be a compact semitopological semigroup with involution, then $C(P)^{**}$ is a quantum semigroup with involution. The coinvolution is of course defined as $(Sf)(t)=f(t^*)$ for $f\in D(S)=C(P)$, $t\in P$. It is easily seen that $S$ satisfied conditions (1)--(4).
\end{example}

On the class of QSI, there exists a duality map ${\cal D}: M\mapsto \hat M$. If $G$ is a locally compact group, then $C_0(G)^{**}$ and $C^*(G)^{**}$ are dual one to another in this sense. For LCQG, there is an analogous duality between {\it universal} algebras.

On the other hand, the duality of QSI connects the two compactifications considered above, at least in the classical case. Let $G$ be a locally compact group and let $\cal D$ be the duality of QSI. Then
$$
{\cal D}^2\big(\,C(G^{WAP})^{**}\big) = C(\b G)^{**}.
$$

\end{document}